# Comments on "Routh Stability Criterion"

T.D.Roopamala  
Assistant Professor  
Sri Jayachamarajendra college of Engineering
.

S.K.Katti  
Professor  
Sri Jayachamarajendra college of Engineering
.

*Abstract—* **In this note, we have shown special case on Routh stability criterion, which is not discussed, in previous literature. This idea can be useful in computer science applications.**

*Keywords- Routh stability criterion, Routh array, Hurwitz criterion, stability.*

I. INTRODUCTION

The Routh stability criterion [1] is an analytical procedure for determining if all the roots of a polynomial have negative real parts, and it is used in the stability analysis of linear time-invariants systems [6]. This stability criterion is useful in various engineering applications [7-8]. There are various special cases discussed in various literature based on Routh criterion [2, 5]. In this small note, one special case based on Routh criterion is considered which is not dealt by the previous authors, to the best of our knowledge.

II. ROUTH STABILITY CRITERION

In order to ascertain the stability of a linear time-invariant system, it is necessary to determine if any of the roots of its characteristics equation lie in the right half of the s-plane. A. Hurwitz and E. J. Routh independently published the method of investigating the sufficient conditions of stability of a system [1]. The Hurwitz criterion is in terms of determinants and the Routh criterion is in terms of array formulation. A necessary and sufficient condition for stability is that all of the elements in the first column of the Routh array be positive. If this condition is not met, the system is unstable, and the number of sign changes in the elements of the first column of the Routh array corresponds to the number of roots of the characteristics equation in the right half of the s-plane. However, the standard procedure fails if we encounter any of the following situations in the formulation of the array [4].

[1] A row of all zeros appears
[2] First element of a row, appearing in first column of the array is zero, but the entire row is not all zeros.

III. PROPOSED SPECIAL CASE ON ROUTH CRITERION

Consider the following polynomial

$$\lambda^4 + 1 = 0 \qquad (1)$$

Applying Routh criterion to above polynomial, we get

| $\lambda^4$ | 1 | 0 | 1 |
|---|---|---|---|
| $\lambda^3$ | $\in$ | $\in$ | $\in$ |
| $\lambda^2$ | -1 | 0 | 0 |
| $\lambda^1$ | $\in$ | $\in$ | 0 |
| $\lambda^0$ | 1 | 0 | 0 |

(2)

In above problem, first row of the Routh's array does not posses all non-zero elements, and the immediate second row has all the elements zeros and hence the problem considered in this note is different from the other existing cases. Now Replacing all the elements in the second row by $\in$, as shown in (2), where $\in$ is a small positive number. Then, we apply Routh method to remaining array formulation. So, in the first column of the Routh array, there exist two changes in sign and hence, this polynomial has complex conjugate eigenvalues with positive real parts.

The actual roots of a polynomial are

$$\begin{aligned}\lambda_{1,2} &= -0.7171 \pm 0.7071j, \\ \lambda_{3,4} &= 0.7171 \pm 0.7071j\end{aligned} \qquad (3)$$

ACKNOWLEDGMENT

We are thankful to Dr. Yogesh V. Hote for suggestion in writing this paper.